\documentclass[preprint,12pt]{elsarticle}
\usepackage{amssymb}
\usepackage{amsfonts}
\usepackage{amsmath}

\newtheorem{theorem}{Theorem}

\newtheorem{corollary}[theorem]{Corollary}
\newtheorem{lemma}[theorem]{Lemma}
\newtheorem{remark}[theorem]{Remark}
\newenvironment{proof}[1][Proof]{\noindent\textbf{#1.} }{\ \rule{0.5em}{0.5em}}

\journal{Statistics \& Probability Letters}

\begin{document}

\begin{frontmatter}

\title{Local asymptotics for the time
of first return to the origin of transient random walk}
\author{R. A. Doney}
\ead{ron.doney@manchester.ac.uk}
\address{University of Manchester}
\author{D. A. Korshunov}
\ead{korshunov@math.nsc.ru}
\address{Sobolev Institute of Mathematics}

\begin{abstract}
We consider a transient random walk on ${\mathbb{Z}}^d$ which is
asymptotically stable, without centering, in a sense which allows different
norming for each component. The paper is devoted to the asymptotics of the
probability of the first return to the origin of such a random walk at
time $n$.
\end{abstract}


\end{frontmatter}

% \linenumbers

\section{Introduction}

Let $S_n$, $n\ge0$, be a random walk in ${\mathbb{Z}}^d$ generated by
independent identically distributed steps $\xi_n=(\xi_{n1},\ldots,\xi_{nd})$%
, $n\ge1$, that is, $S_0=0$, $S_n=\xi_1+\ldots+\xi_n$. Denote $\tau_0=0$ and
recursively $\tau_{n+1}=\min\{k>\tau_n:S_k=0\}$; by standard convention $%
\min\emptyset=\infty$. Then $\tau=\tau_1$ is the first return to the origin
of the random walk $S_n$.

In this paper we study the asymptotic behaviour of
$$
p_{n}:={\mathbb{P}}\{\tau =n\}
={\mathbb{P}}\{\mbox{the first return to zero
occurs at time }n\}
$$
as $n\rightarrow \infty $. Put $G(B)={\mathbb{P}}\{\tau \in B\}$ and
$$
p:=\sum_{n=1}^{\infty }p_{n}={\mathbb{P}}\{\tau <\infty \}=G[1,\infty )\leq
1.
$$%
The measure $u_{n}:={\mathbb{P}}\{S_{n}=0\}$ on ${\mathbb{Z}}^{+}$ is
actually the renewal measure generated by the $\tau $'s, that is,
\begin{equation}
u_{n}=\sum_{k=0}^{\infty }{\mathbb{P}}\{\tau _{k}=n\}=\sum_{k=0}^{\infty
}G^{\ast k}(n).  \label{ren.u}
\end{equation}%
Then
\begin{equation}
\sum_{n=0}^{\infty }{\mathbb{P}}\{S_{n}=0\}=\sum_{k=0}^{\infty }G^{\ast
k}[0,\infty )=\sum_{k=0}^{\infty }p^{k}=\frac{1}{1-p},
\end{equation}%
which implies
\begin{equation}
p=\frac{\sum_{n=1}^{\infty }{\mathbb{P}}\{S_{n}=0\}}{1+\sum_{n=1}^{\infty }{%
\mathbb{P}}\{S_{n}=0\}}.
\end{equation}

The random walk $S_n$ is called aperiodic if $Z^d$ is a minimal lattice for $%
S_n$ in the sense that, for every $\varepsilon>0$,
$$
\sup_{\lambda\in[-\pi,\pi]^d\setminus[-\varepsilon,\varepsilon]^d} |{\mathbb{%
E}}e^{i(\lambda,\xi_1)}| < 1.
$$
Aperiodicity is clearly no essential restriction, as the state space can
always be redefined, if necessary, so as to make a random walk aperiodic.

In the sequel we will be studying aperiodic random walks on $\mathbb{Z}^d$
which are \emph{asymptotically stable} in the following sense: there is
sequence $c_n=(c_{n1},\ldots,c_{nd})$ such that
$$
X_n:=(S_{n1}/c_{n1},\ldots,S_{nd}/c_{nd})\overset{D}{\to }Y=(Y_1,\ldots,Y_d),
$$
where $Y$ is a strictly $d$-dimensional stable random variable. Since this
implies that each component of $X_n$ is asymptotically stable, we know that
each $c_{nr}$ is in the class $RV(1/\alpha_r)$
of regularly varying at infinity with index $1/\alpha_r$ sequences
(see, e.g. \cite[Section 1.9]{BGT}), where $\alpha_r\in (0,2]$ is the
index of the univariate stable random variable $Y_r$. Thus $C_n:=\prod_1^d
c_{nr}$ is in $RV(\eta)$, where $\eta=\sum_1^d 1/\alpha_r\ge d/2$. We need
the following local limit theorem, in which $g$ denotes the density function
of $Y$.

\begin{theorem}
If $S_n$ is an aperiodic random walk on $\mathbb{Z}^d$ which is
asymptotically stable in the above sense, it holds that uniformly for $x\in%
\mathbb{Z}^d$
$$
C_n\mathbb{P}\{S_n=x\}=g(x_1/c_{n1},\ldots,x_d/c_{nd})+o(1) \text{ as }%
n\to\infty.
$$
In particular, $u_n=\mathbb{P}\{S_n=0\}\sim g(0,\ldots,0)/C_n$
as $n\to\infty$.
\end{theorem}

For $d=1$ this is the classical local limit theorem of Gnedenko, see \cite[%
\S ~50]{KG}; for the case $d=2$ it is proved in \cite{Doney}, and as
remarked there, the proof extends in a straightforward way
to the case $d>2$. It can also be viewed as a special case of
Theorem 6.4 in \cite{Griffin},
where the more general case of matrix norming is treated.

Since it is known that $g$ and its derivatives are bounded,
we deduce the following

\begin{corollary}\label{C}
If $S_n$ is an aperiodic random walk on $\mathbb{Z}^d$ which is
asymptotically stable, there exists a constant $k$ such that
\begin{eqnarray*}
{\mathbb{P}}\{S_n=x\} &\le& k/C_n\quad \mbox{ for all }x\in {\mathbb Z}^d,
\end{eqnarray*}
and, for every fixed $k$,
$$
{\mathbb{P}}\{S_{n-k}=x\}={\mathbb{P}}\{S_n=x\}+o(1/C_n)
$$
as $n\to\infty$ uniformly for $x\in\mathbb{Z}^d$.
\end{corollary}

\begin{remark}
Since every random walk is transient when $d\ge 3$, a result which for the
simplest symmetric random walk goes back to Polya \cite{Polya}, the
requirement of transience only features for $d=1$ and $d=2$, when under our
assumptions it is equivalent to $\sum_{n=1}^\infty 1/C_n<\infty $.
\end{remark}

Suppose we know that ${\mathbb{P}}\{\tau=n\}\in RV(-\gamma)$ for some $%
\gamma\ge 1$. Then $G/p$ is so-called \textit{locally subexponential}
distribution. In this case, the local asymptotics of the defective renewal
function is described in \cite[Proposition 12]{AFK}: the following limit
exists:
$$
\lim_{n\to\infty}\frac{\sum_{k=0}^\infty G^{\ast k}(n)}{G(n)}%
=\sum_{k=0}^\infty kp^{k-1}=\frac{1}{(1-p)^2},
$$
so that ${\mathbb{P}}\{S_n=0\}\sim {\mathbb{P}}\{\tau =n\}/(1-p)^2$. This is
the intuition behind the following, which is our main result.

\begin{theorem}\label{d.ge.3}
Let $S_n$ be an aperiodic, transient random walk on $\mathbb{Z%
}^d$, $d\ge 1$, which is asymptotically stable in the above sense. Then as $%
n\to\infty $,
$$
{\mathbb{P}}\{\tau=n\}\sim (1-p)^2{\mathbb{P}}\{S_n=0\} \sim \frac{%
(1-p)^2g(0,\ldots,0)}{C_n}.
$$
\end{theorem}

\begin{remark}
The special case where $d\ge 3,{\mathbb{E}}\xi _{n}=0$ and ${\mathbb{E}}\xi
_{nj}\xi _{nj}=B_{ij}<\infty $ with $\det B\not=0$ reads
$$
{\mathbb{P}}\{\tau =n\}\sim (1-p)^{2}{\mathbb{P}}\{S_{n}=0\}\sim \frac{%
(1-p)^{2}}{(2\pi )^{d/2}\sqrt{\det B}}n^{-d/2}.
$$%
This was proved in an unpublished communication by one of us, and is quoted
in Chapter A.6 of \cite{G}; this illustrates the increasing importance of
local results such as this in Mathematical Physics.
To the best of our knowledge, Theorem \ref{d.ge.3} was not proved
in the literature even in the case of the simplest symmetric
random walk on ${\mathbb Z}^d$ with $d\ge3$.
\end{remark}

In the next section we prove Theorem \ref{d.ge.3} by analytic means via a
Banach algebra technique; this is the method that was used in proving the
above special case. Then in Section \ref{sec:prob} we give a probabilistic
proof capturing the most probable way that large values of $\tau $ occur.

Note also that the recurrent one-dimensional case $d=1$ was first studied by
Kesten in \cite{Kesten} where it was proved in Theorems 7 and 8 that when $%
\mathbb{E\xi =}0$ and $\mathbb{E\xi }^2<\infty$, we have ${\mathbb{P}}%
\{\tau=n\}\sim \sqrt{\frac{\mathrm{Var}\xi}{2\pi }}n^{-3/2}$ as $n\to\infty$%
; the recurrent case of convergence to a stable law with index $%
\alpha\in\lbrack 1,2]$ was also considered. Two different approaches for
proving this equivalence may be found in \cite[Theorem 1.2]{BLPW}.

In dimension 2 Jain and Pruitt \cite[Theorem 4.1]{JP} proved, assuming zero
mean and finite covariance $B$, that ${\mathbb{P}}\{\tau =n\}\sim 2\pi \sqrt{%
\det B}n^{-1}\log ^{-2}n$ as $n\to\infty$.

The only local result for
dimensions 3 and higher we found is one by Kesten and Spitzer \cite[Theorem
1b]{KestenSpitzer} where they proved that ${\mathbb{P}}\{\tau=n+1\}\sim {%
\mathbb{P}}\{\tau=n\}$ as $n\to\infty$ given ${\mathbb{E}}\xi=0$.

In conclusion note that different mechanisms are involved in formation of
large deviations of $\tau$ in dimensions $d=1$, $d=2$ and $d\ge 3$. In the
one dimensional recurrent case, the equation \eqref{ren.u} says that we deal
with the renewal process generated by $\tau$'s where $\tau$ has an infinite
mean. In principle the same is true for the case $d=2$ but here the tail of $%
\tau$ is very heavy, it is slowly varying at infinity. In the case $d\ge 3$
transience holds, so that we have the renewal process generated by a
defective distribution which yields that the renewal atoms $u_n$ are
proportional to $p_n$; this is the main topic addressed to in the present
article.

\section{Banach algebra approach}

The proof of Theorem \ref{d.ge.3} is straightforward, if we assume
additionally that $\sum_{n=1}^\infty u_n<1$ which is equivalent to $p<1/2$.
The discrete renewal relation \eqref{ren.u} between the $u$'s and $p$'s
implies that
$$
1+u(s)=\frac{1}{1-p(s)},
$$
where we put $u(s)=\sum_{n=1}^\infty u_ns^n$ and $p(s)=\sum_{n=0}^\infty
p_ns^n$, so that
\begin{eqnarray*}
p(s) &=& \frac{u(s)}{1+u(s)} = \sum_{n=1}^\infty (-1)^{n+1}(u(s))^n,\quad
|s|\le1.
\end{eqnarray*}
This is equivalent to
$$
p_n=\sum_{k=1}^n (-1)^{k+1} u_n^{*(k)}
$$
where $u_n^{*(k)}$ is the $k$-fold convolution of $\{u_n\}$. Note that $u_n$
is regularly varying at infinity so that the defective distribution $\{u_n\}$
is locally subexponential, see \cite{AFK} or \cite[Chapter 4]{FKZ}. Then, as
follows from \cite[Theorem 4.30]{FKZ},
$$
\frac{p_n}{u_n} = \sum_{k=1}^n (-1)^{k+1}\frac{u_n^{*(k)}}{u_n} \to
\sum_{k=1}^\infty (-1)^{k+1}k(u(1))^{k-1} =\frac{1}{(1+u(1))^2}.
$$
Note $1+u(1)=1/(1-p)$, so $p_n\sim (1-p)^2u_n$ as $n\to\infty$.

What happens if $\sum_{n=1}^{\infty }u_n\ge 1$? In this case we can't
expand $u(s)/(1+u(s))$ as a power series in $u(s)$. Nevertheless this is an
analytic function in $u(s)$, for all (complex) $s$ in $|s|\le 1$, because $%
u(1)<\infty $ by the assumed transience. We can then apply Theorem 1 from
the paper \cite{CNW} and get the same result; this reference is based on
Banach algebra techniques.

\section{Probabilistic approach}

\label{sec:prob}

The starting point of the probabilistic proof of Theorem \ref{d.ge.3} is the
following result, which holds in any dimension.

\begin{lemma}
\label{lemma1} Let a sequence $Q_n$ which is regularly varying at infinity
be such that
\begin{eqnarray}
{\mathbb{P}}\{S_n=x\} &\le &Q_n\quad \mbox{ for all }x\in {\mathbb{Z}}^d,
\label{Qn.Sn} \\
\sum_{n=1}^{\infty }Q_n &<&\infty ,  \label{sum Qn}
\end{eqnarray}%
and, for every fixed $k$,
\begin{equation}
{\mathbb{P}}\{S_{n-k}=x\}={\mathbb{P}}\{S_n=x\}+o(Q_n)  \label{n-k.n}
\end{equation}%
as $n\to \infty $ uniformly in all $x\in {\mathbb{Z}}^d$.
Let $r_n$ be any fixed unboundedly increasing sequence. Then
$$
{\mathbb{P}}\{S_n=x,\tau >n\}=(1-p){\mathbb{P}}\{S_n=x\}+o(Q_n)
$$%
as $n\to \infty $ uniformly in $x\in {\mathbb{Z}}^d$ such that $\Vert x\Vert
>r_n$.
\end{lemma}

\proof It is equivalent to prove the relation
\begin{equation}
{\mathbb{P}}\{S_n=x,\tau <n\}=p{\mathbb{P}}\{S_n=x\}+o(Q_n).
\label{lemma1.eq}
\end{equation}%
We start with the following decomposition:
\begin{eqnarray*}
{\mathbb{P}}\{S_n=x,\tau <n\} &=&\sum_{k=1}^{n-1}{\mathbb{P}}\{S_n=x,\tau
=k\} \\
&=&\sum_{k=1}^{n-1}{\mathbb{P}}\{S_{n-k}=x\}{\mathbb{P}}\{\tau =k\}.
\end{eqnarray*}%
For every fixed $N$, we have
\begin{eqnarray}
\sum_{k=N}^{n-N}{\mathbb{P}}\{S_{n-k}=x\}{\mathbb{P}}\{\tau =k\} &\le
&\sum_{k=N}^{n-N}{\mathbb{P}}\{S_{n-k}=x\}{\mathbb{P}}\{S_{k}=0\}  \notag \\
&\le &\sum_{k=N}^{n-N}Q_{k}Q_{n-k},  \label{est1}
\end{eqnarray}%
by the condition \eqref{Qn.Sn}. For every fixed $k$, ${\mathbb{P}}%
\{S_{k}=x\}\to 0$ as $\Vert x\Vert \to \infty $. Together with \eqref{Qn.Sn}
it implies that, for every fixed $N$,
\begin{eqnarray}
\sum_{k=n-N}^{n-1}{\mathbb{P}}\{S_{n-k}=x\}{\mathbb{P}}\{\tau =k\}
&\le& \sum_{k=n-N}^{n-1}{\mathbb{P}}\{S_{n-k}=x\}Q_{k}\nonumber\\
&\le& \sum_{k=n-N}^{n-1}Q_{k}o(1)=o(Q_n)  \label{est2}
\end{eqnarray}
as $n\to \infty $ uniformly in $\Vert x\Vert \ge r_n$. Finally, by %
\eqref{n-k.n}, for every fixed $N$,
\begin{equation}
\sum_{k=1}^{N}{\mathbb{P}}\{S_{n-k}=x\}{\mathbb{P}}\{\tau =k\}={\mathbb{P}}
\{S_n=x\}\sum_{k=1}^{N}{\mathbb{P}}\{\tau =k\}+o(Q_n)  \label{est3}
\end{equation}%
as $n\to \infty $ uniformly in all $x$.

Combining \eqref{est1}--\eqref{est3} we obtain that
\begin{eqnarray*}
\lefteqn{\limsup_{n\to \infty }\frac{1}{Q_n}|\sum_{k=1}^{n-1}\mathbb{P}%
\{S_n=x,\tau =k\}-\mathbb{P}\{S_n=x\}\mathbb{P}\{\tau \le N\}|} \\
&&\hspace{70mm}\le \limsup_{n\to \infty }\frac{1}{Q_n}%
\sum_{k=N}^{n-N}Q_{k}Q_{n-k}.
\end{eqnarray*}%
By regular variation of $Q_n$ at infinity, taking into account that the
series \eqref{sum Qn} converges, we conclude that the right hand side can be
made as small as we please by choosing $N$ sufficiently large; this is also
a well-known property in the theory of locally subexponential distributions,
see e.g \cite[Chapter 4]{FKZ}. Now the proof of \eqref{lemma1.eq} follows by
letting $N\to \infty $.

\begin{lemma}
\label{lemma2} Under the conditions of Lemma \ref{lemma1},
as $n\to \infty $,
$$
{\mathbb{P}}\{\tau =n\}=(1-p)^{2}{\mathbb{P}}\{S_n=0\}+o(Q_n).
$$
\end{lemma}

\proof Let $m$ be such that $n=2m$ in the case of even $n$ and $n=2m+1$
otherwise. We have
$$
{\mathbb{P}}\{\tau =n\}=\sum_{x\in {\mathbb{Z}}^d\setminus 0}{\mathbb{P}}%
\{S_m=x,\tau =n\}.
$$%
Since $S_m$ and the sequence $\{S_{m+k}-S_m,k\ge 1\}$ are independent, the $%
x$th summand in the latter sum is equal to the product
\begin{eqnarray*}
\lefteqn{\mathbb{P}\{S_m=x,\tau >m\}} \\
&&\times {\mathbb{P}}\{S_{m+k}-S_m\not=-x\text{ for all }k=1,\ldots
,n-m-1,S_n-S_m=-x\}.
\end{eqnarray*}%
The second probability here is equal to
\begin{eqnarray*}
\lefteqn{\mathbb{P}\{S_n-S_{m+k}\not=0\text{ for all }k=1,\ldots
,n-m-1,S_n-S_m=-x\}} \\
&=&{\mathbb{P}}\{\widetilde{S}_{k}\not=0\text{ for all }k=1,\ldots ,n-m-1,%
\widetilde{S}_{n-m}=-x\} \\
&=&{\mathbb{P}}\{\widetilde{\tau }>n-m,\widetilde{S}_{n-m}=-x\},
\end{eqnarray*}%
where $\widetilde{S}_{k}:=\widetilde{\xi }_{1}+\ldots +\widetilde{\xi }_{k}$%
, $\widetilde{\xi }_{k}:=\xi _{n-k+1}$, and $\widetilde{\tau }$ is the first
return time to zero of the random walk $\widetilde{S}_{k}$. The random walk $%
\widetilde{S}_{k}$ has the same distribution as $S_{k}$, so
\begin{eqnarray}
{\mathbb{P}}\{\tau =n\} &=&\sum_{x\in {\mathbb{Z}}^d\setminus 0}{\mathbb{P}}%
\{S_m=x,\tau >m\}{\mathbb{P}}\{\widetilde{S}_{n-m}=-x,\widetilde{\tau }>n-m\}
\notag  \label{sigmas} \\
&=&\sum_{x\in {\mathbb{Z}}^d\setminus 0}{\mathbb{P}}\{S_m=x,\tau >m\}{%
\mathbb{P}}\{S_{n-m}=-x,\tau >n-m\}  \notag \\
&=&\Sigma _{1}+\Sigma _{2},
\end{eqnarray}
where $\Sigma_1$ is the sum over $\Vert x\Vert \le \log n=:r_n$ and $\Sigma_2
$ is the sum over $\Vert x\Vert >\log n$.
By the condition \eqref{Qn.Sn} and regular variation of $Q_n$,
\begin{eqnarray}\label{sigma1}
\Sigma _{1} &\le&
\sum_{\Vert x\Vert \le \log n}{\mathbb{P}}\{S_m=x\}
{\mathbb{P}}\{S_{n-m}=-x\}\nonumber\\
&\le& Q_mQ_{n-m}(2\log n)^d=o(Q_n)
\quad \mbox{ as }n\to \infty.
\end{eqnarray}
By Lemma \ref{lemma1}, as $n\to\infty$,
\begin{eqnarray*}
\Sigma _{2} &=&(1-p)\sum_{\Vert x\Vert >\log n}{\mathbb{P}}\{S_m=x,\tau >m\}[%
{\mathbb{P}}\{S_{n-m}=-x\}+o(Q_n)] \\
&=& (1-p)\sum_{\Vert x\Vert >\log n}{\mathbb{P}}\{S_m=x,\tau >m\}{\mathbb{P}}%
\{S_{n-m}=-x\}+o(Q_n).
\end{eqnarray*}
Repeating these arguments to the first multiple we obtain that
$$
\Sigma _{2}=(1-p)^{2}\sum_{\Vert x\Vert >\log n}{\mathbb{P}}\{S_m=x\}{%
\mathbb{P}}\{S_{n-m}=-x\}+o(Q_n).
$$%
Taking also into account \eqref{sigma1} we finally obtain that
\begin{equation}
\Sigma _{2}=(1-p)^{2}{\mathbb{P}}\{S_n=0\}+o(Q_n).  \label{14}
\end{equation}%
Substituting \eqref{sigma1} and (\ref{14}) into \eqref{sigmas}, we arrive at
the desired conclusion.

The proof of Theorem \ref{d.ge.3} is now immediate, since by Corollary
\ref{C} the sequence $Q_n=k/C_n$ satisfies all conditions
of Lemma \ref{lemma1}, for suitable $k$. The proof is complete.

\begin{remark}
The same question may be addressed in the more general setting of matrix
norming with the help of results by Griffin \cite{Griffin}. The key point is
that the norming sequence $|B_n|$ in his Theorem 6.4 is automatically
regularly varying: we owe this comment to Phil Griffin, in a private
communication. It is then easy to see that our result extends to this
situation whenever $d\ge 3$ or $d=2$ and transience is assumed.
\end{remark}

\section*{Acknowledgment}

The authors are grateful to the referees for their
careful reading and for helpful comments and remarks.


\begin{thebibliography}{99}
\bibitem{AFK}
Asmussen, S., Foss, S., Korshunov, D. (2003).
Asymptotics for sums of random variables with local subexponential behaviour.
\textit{J. Theoret. Probab.} \textbf{16} 489--518.

\bibitem{BLPW}
Bender, E. A., Lawler, G. F., Pemantle, R., Wilf, H. S. (2004).
Irreducible compositions and the first return to the origin of a
random walk.
\textit{S\'eminaire Lotharingien de Combinatoire} \textbf{50}
Article B50h.

\bibitem{BGT}
Bingham, N. H., Goldie, C. M., and Teugels J. L. (1987).
\textit{Regular variation},
Cambridge University Press, Cambridge.

\bibitem{CNW}
Chover, J., Ney, P., and Wainger, S. (1973).
Functions of probability measures.
\textit{J. d'Analyse Math\'ematique} \textbf{26} 255--302.

\bibitem{Doney}
Doney, R. (1991).
A bivariate local limit theorem.
\textit{J. Multivariate Anal.} \textbf{36} 95--102.

\bibitem{FKZ}
Foss, S., Korshunov, D., Zachary, S. (2011).
\textit{An Introduction to Heavy-Tailed and Subexponential Distributions}. Springer,
New York, to appear.

\bibitem{G}
Giacomin, G. (2007).
\textit{Random Polymer Models}.
Imperial College Press, London.

\bibitem{Griffin}
Griffin, P. S. (1986).
Matrix normalized sums of independent
identically distributed random variables.
\textit{Ann. Probab.} \textbf{14} 224--246.

\bibitem{JP}
Jain, N. C. and Pruitt, W. E. (1972).
The range of random walk.
\textit{Proc. Sixth Berkeley Symp. Math. Statist. Probab.} \textbf{3}
31---50. Univ. California Press, Berkeley.

\bibitem{Kesten}
Kesten, H. (1963).
Ratio theorems for random walks II.
\textit{J. d'Analyse Math\'ematique} \textbf{9} 323--379.

\bibitem{KestenSpitzer}
Kesten, H., Spitzer, F. (1963).
Ratio theorems for random walks I.
\textit{J. d'Analyse Math\'ematique} \textbf{9} 285--322.

\bibitem{KG}
Kolmogorov, A. N., Gnedenko, B. V. (1954).
\textit{Limit distributions for sums of independent random variables}.
Addison-Wesley Publishing Company, Reading.

\bibitem{Polya}
Polya, G. (1921).
\"{U}ber eine Aufgabe der Wahrscheinlichkeitsrechnung
betreffend die Irrfahrt im Strassennetz.
\textit{Mathematische Annalen} \textbf{84} 149--160.
\end{thebibliography}
\end{document}